\documentclass[12pt]{article}
\usepackage{amssymb,amsmath}

\renewcommand{\theequation}{\thesection.\arabic{equation}}
\numberwithin{equation}{section}

\input amssym.def

\newcommand{\totimes}{{\widetilde{\otimes}}}
\newcommand{\bul}{{\bullet}}
\newcommand{\al}{{\alpha}}
\newcommand{\la}{{\lambda}}

\newcommand{\h}{{\hbar}}

\newcommand{\mP}{{\mathfrak{P}}}

\newcommand{\mb}{{\mathfrak{b}}}

\newcommand{\md}{{\mathfrak{d}}}

\newcommand{\mK}{{\mathfrak{K}}}

\newcommand{\mE}{{\mathfrak{E}}}

\newcommand{\mD}{{\mathfrak{D}}}
\newcommand{\mQ}{{\mathfrak{Q}}}

\newcommand{\by}{{\tilde y}}
\newcommand{\bz}{{\tilde z}}

\newcommand{\om}{{\omega}}
\newcommand{\Om}{{\Omega}}

\newcommand{\si}{{\sigma}}

\newcommand{\ka}{{\kappa}}
\newcommand{\G}{{\Gamma}}
\newcommand{\pa}{{\partial}}

\newcommand{\W}{{\cal W}}

\newcommand{\cC}{{\cal C}}

\newcommand{\cD}{{\cal D}}

\newcommand{\cH}{{\cal H}}

\newcommand{\cO}{{\cal O}}

\newcommand{\bbA}{{\Bbb A}}

\newcommand{\bbC}{{\Bbb C}}
\newcommand{\bbR}{{\Bbb R}}
\newcommand{\bbZ}{{\Bbb Z}}

\newcommand{\n}{{\nabla}}
\newcommand{\te}{\theta}

\newcommand{\de}{{\delta}}

\newcommand{\OmW}{\Om(M,\W)}

\newcommand{\OmC}{\Om(M,\cC)}

\date{}
\newtheorem{defi}{Definition}
\newtheorem{pred}{Proposition}

\newtheorem{teo}{Theorem}
\newtheorem{cor}{Corollary}
\begin{document}

\vspace{-10cm}
\begin{flushright}
 \begin{minipage}{1.2in}
 ITEP-TH-15/04
 \end{minipage}
\end{flushright}
\vspace{1.5cm}

\begin{center}
{\Large\bf Hochschild Cohomology versus \\[0.3cm] De Rham Cohomology
\\[0.3cm] without Formality Theorems.}\\[1cm]
V.A. Dolgushev\footnote{On leave of
absence from: ITEP (Moscow)} \\[0.3cm]
{\it Department of Mathematics, MIT,} \\
{\it 77 Massachusetts Avenue,} \\
{\it Cambridge, MA, USA 02139-4307,}\\
{\it E-mail address: vald@math.mit.edu}
\end{center}

\begin{abstract}
We exploit the Fedosov-Weinstein-Xu (FWX) resolution
proposed in q-alg/9709043
to establish an isomorphism between the ring of
Hochschild cohomology of the quantum algebra of functions
on a symplectic manifold $M$ and the ring $H^{\bul}(M,\bbC((\h)))$ of
De Rham cohomology of $M$ with the coefficient field $\bbC((\h))$
without making use of any version of the formality theorem.
We also show that the Gerstenhaber bracket induced on $H^{\bul}(M,\bbC((\h)))$
via the isomorphism is vanishing.
We discuss equivariant properties of the isomorphism
and propose an analogue of this statement in
an algebraic geometric setting.
\end{abstract}

\section{Introduction}
Although Kontsevich's formality theorem \cite{K},\cite{K1}, \cite{Dima}
and its generalizations \cite{B-Gilles},\cite{Damien}, \cite{CFT}, \cite{CEFT}, \cite{FTC},
\cite{F-Sh}, \cite{Sh}, \cite{Sh1}, \cite{TT}, \cite{TT1}
give immediate solutions of various problems of deformation quantization,
the complicated technique of homotopy algebraic structures is not required
for many of these problems.
In fact, a number of examples show that for many problems of deformation quantization
this technique can be replaced by simpler
arguments \cite{B-K}, \cite{Wick},
\cite{Fedosov}, \cite{Fedosov1}, \cite{Vanya}, \cite{L-S}, \cite{NT}, \cite{NT1},
\cite{RT}, \cite{WX}, \cite{Xu}.

In paper \cite{WX} by A. Weinstein and P. Xu the authors proposed
a resolution of the vector space of local Hochschild
cochains of the ring of quantum functions on a symplectic manifold $M$.
They used this resolution to
prove that the graded vector space of Hochschild cohomology of the ring
of quantum functions on a symplectic manifold is isomorphic to
the graded vector space $H^{\bul}(M,\bbC((\h)))$ of De Rham cohomology
of $M$ with the coefficients in the field $\bbC((\h))$
of formal Laurent power series.
In our paper we use this resolution
to show that the above cohomology spaces are isomorphic as
rings. In principle, the compatibility
of Kontsevich's quasi-isomorphism \cite{K} with the
cup-product in Hochschild cohomology implies
an analogous statement for a general Poisson
manifold. However, it is instructive to see that in
the symplectic case, the result can be proven by a
much simpler argument.
From our considerations it will also be transparent that
the Gerstenhaber bracket induced
on De Rham cohomology via the above isomorphism
with Hochschild cohomology is vanishing.

A dual version of this assertion was proven in paper
\cite{NT} by Nest and Tsygan on the algebraic index theorem.
Namely, in \cite{NT} the authors describe
Hochschild and cyclic homology of the algebra of quantum
functions with compact support
on an arbitrary symplectic manifold using
the spectral sequence associated to the
$\h$-adic filtration.

In this context, it is worth mentioning
paper \cite{B-G}, in which similar results were
obtained for Hochschild and cyclic homology
of the ring of pseudo-differential symbols on an
arbitrary smooth manifold.

The structure of this paper is as follows.
In the next section we formulate the main results of
the paper (See theorems $1$ and $2$).
In the third section we remind the Fedosov
deformation quantization, recall the Fe\-do\-sov-Wein\-stein-Xu
resolution of the algebra of local Hochschild cochains
of the quantum ring of functions on a symplectic manifold,
and prove theorems $1$ and $2$ formulated in the previous
section. In the concluding section we discuss applications and
the variations of theorems $1$ and $2$\,. In Appendices A and B
at the end of the paper we recall algebraic structures on
Hochschild cochains and propose an equivariant
homotopy operator for the cohomological Hochschild
complex of the formal Weyl algebra.

Throughout the paper we assume the summation over repeated indices.
 We omit symbol $\wedge$ referring to a local basis
of exterior forms as if we thought of $dx^i$'s as anti-commuting variables.
We denote by $d$ the De Rham differential. By a nilpotent linear operator
we always mean an operator whose second power is vanishing.

\section{Algebra of local Hochschild cochains of the
ring of quantum functions on a symplectic manifold.}
Let $M$ be an even dimensional smooth manifold endowed with a symplectic form
$\om=\om_{ij}(x)dx^i\wedge dx^j$. Here $i,j$ run from $1$ to $2n=dim M$
and $x^i$ denote local coordinates. Let
$\om^{ij}(x)\pa_{x^i}\wedge \pa_{x^j}$ denote the respective
Poisson tensor
$$
\om^{ik}(x)\om_{kj}(x)=\de^i_j\,.
$$

According to the standard terminology of deformation
quantization \cite{Bayen}, \cite{Ber} a star-product
on the symplectic manifold $M$ is an associative
non-commutative $\bbC((\h))$-linear product in
$C^{\infty}(M)((\h))$ given by the formal power
series of bidifferential operators $S_k$
\begin{equation}
\label{Berezin}
a\ast b = ab + \sum_{k=1}^{\infty}\h^k S_{k}(a,b)\,,
\qquad a,b\in C^{\infty}(M)
\end{equation}
and such that
$$
a\ast b - b \ast a = \h \{a,b\} ~ mod~ \h^2\,, \qquad
a\ast 1 = 1 \ast a = a.
$$
where $\{,\}$ denotes the Poisson bracket associated to
the symplectic structure on $M$\,.

Furthermore, two star-products $\ast$ and $\ast'$ are
called equivalent if there exists a formal series
$$
Q = I + \sum_{k=1}^{\infty}\h^k Q_k
$$
of differential operators $Q_k$ on $M$ such that

\begin{equation}
\label{P-equiv}
Q( a \ast b ) = (Q a)\ast' (Q b)\,,
\qquad \forall~ a,b \in C^{\infty}(M)((\h))\,.
\end{equation}

We denote by $\bbA$ the algebra $C^{\infty}(M)((\h))$ of formal
Laurent power series of complex valued smooth functions
on $M$ with the multiplication (\ref{Berezin}).

By definition the vector space $C_{loc}^{k}(\bbA)$
of {\it local Hochschild
$k$-cochains} of the algebra $\bbA$ is the subspace of
$C_{loc}^{k}(\bbA)\subset Hom_{\bbC((\h))}(\bbA^{\otimes\,k}, \bbA)$ of
$\bbC((\h))$-linear homomorphisms from $\bbA^{\otimes\,k}$ to
$\bbA$ preserving supports of functions. In
\cite{Cahen} it is shown that any polylinear map
$L\in  Hom_{\bbC}((C^{\infty}(M))^{\otimes\,k}, C^{\infty}(M))$
which preserves supports of functions
can be locally represented as a polydifferential
operator\footnote{If the manifold $M$ is compact
the word ``locally'' can be omitted.}.
Thus any element $P\in C_{loc}^k(\bbA)$ can be
defined locally as the following formal series
of $k$-differential operators
\begin{equation}
\label{operr-A}
P(a_1, \dots, a_k) =\sum_{m\in \bbZ}
\h^m P_m(a_1, \dots, a_k)\,,
\end{equation}
where the summation in $m$ is bounded below.

It is not hard to see that the Hochschild differential
$\pa$ (\ref{Hoch-diff}), the cup-product (\ref{cup-W}) and
the Gerstenhaber bracket (\ref{Gerst-W}) preserve
locality of Hochschild cochains and therefore are
well-defined on $C_{loc}^{\bul}(\bbA)$.
In what follows by the ring of Hochschild cohomology
of $\bbA$ we mean the vector space
$$
HH^{\bul}_{loc}(\bbA) = H(C^{\bul}_{loc}(\bbA), \pa)
$$
with the multiplication induced by (\ref{cup-W}).

We have to mention that since $\bbA$ is a deformation
of the algebra of smooth function on $M$ the complex
of all Hochschild cochains of $\bbA$
is a rather intractable object. In particular,
Hochschild cohomology $HH^{k}(A)$ for $k>1$ is
not known even for the ordinary commutative algebra
of smooth functions on $\bbR^n$\,.

Recall that given a symplectic torsion free
connection $\n$ and a Fedosov representative
$\displaystyle \Om_F\in \frac1{\h}\Om^2(M)[[\h]]$
(\ref{repres}) one can construct an algebra
$\bbA_{\n,\Om_F}$ which quantizes the symplectic
manifold $M$ in the sense of (\ref{Berezin}).
It is well-known (see theorem \ref{Xu}) that
for different Fedosov representatives $\Om_F$
the algebras $\bbA_{\n, \Om_F}$ exhaust all
equivalence classes of deformation quantizations
of $M$. In this paper we refer to the triple
$(M, \n, \Om_F)$ consisting of a symplectic manifold
$M$, a symplectic torsion free connection $\n$ and
a Fedosov representative $\Om_F$ (\ref{repres}) as
{\it Fedosov data}.

The main results of the paper are formulated
in the following two theorems
\begin{teo}
For the algebra $\bbA=(C^{\infty}(M)((\h)),\ast)$
of quantum functions on a symplectic manifold $M$
the ring of Hochschild cohomology
$$
HH^{\bul}_{loc}(\bbA)
$$
is isomorphic to the ring of De Rham cohomology
$$
H_{DR}(M)\otimes \bbC((\h))\,,
$$
with coefficients in the field $\bbC((\h))$
of formal Laurent power series.
The Lie bracket induced on
$H^{\bul}_{DR}(M)\otimes \bbC((\h))$
via this isomorphism is vanishing.
\end{teo}

\begin{teo}
To the Fedosov data $(M, \n, \Om_F)$ one can naturally
assign a differential
graded (DG) associative algebra $K^{\bul}_{\n, \Om_F}$ and
a pair of embeddings
\begin{equation}
\label{emb-Rham}
\mE_{\n, \Om_F}\,:\, \Om^{\bul}(M)((\h))
\hookrightarrow K^{\bul}_{\n, \Om_F}\,,
\end{equation}
\begin{equation}
\label{emb-Hoch}
\mD_{\n, \Om_F}\,:\, C^{\bul}_{loc}(\bbA_{\n, \Om_F})
\hookrightarrow K^{\bul}_{\n, \Om_F}\,,
\end{equation}
which are quasi-isomorphisms of the corresponding
DG associative algebras. If $g$ is a diffeomorphism of $M$
then the corresponding embeddings $\mD_{\n, \Om_F}$,
$\mE_{\n, \Om_F}$, $\mD_{g_{\ast}\n, g_{\ast} \Om_F}$,
and $\mE_{g_{\ast}\n, g_{\ast} \Om_F}$ fit into the
commutative diagram
\begin{equation}
\label{al-equiv}
\begin{array}{ccccc}
C^{\bul}_{loc}(\bbA_{\n, \Om_F})  &
\stackrel{\mD_{\n, \Om_F}}{\longrightarrow}&
K^{\bul}_{\n, \Om_F}
& \stackrel{\mE_{\n, \Om_F}}{\longleftarrow}&
\Om^{\bul}(M)\otimes \bbC((\h))\\[0.3cm]
\downarrow^{g_{\ast}}  & ~  & \downarrow^{g_{\ast}}  & ~
& \downarrow^{g_{\ast}} \\[0.3cm]
C^{\bul}_{loc}(\bbA_{g_{\ast}\n,\, g_{\ast}\Om_F})  &
\stackrel{\mD_{g_{\ast}\n,\, g_{\ast} \Om_F}}{\longrightarrow}
&K^{\bul}_{g_{\ast}\n,\, g_{\ast}\Om_F} &
\stackrel{\mE_{g_{\ast}\n,\, g_{\ast} \Om_F}}{\longleftarrow}
&\Om^{\bul}(M)\otimes \bbC((\h))\,.
\end{array}
\end{equation}
\end{teo}
Proofs of the above statements are given in the
next section. The proofs are based on the use
of the Fedosov-Weinstein-Xu (FWX) resolution \cite{WX}
of the vector space of local Hochschild cochains of
$\bbA_{\n, \Om_F}$.
An interesting corollary and a couple of variations of
theorems $1$ and $2$ are discussed in the
concluding section of the paper.

{\bf Remark.} We would like to mention that the triviality of
the operation induced by the Gerstenhaber bracket
on the Hochschild cohomology can be viewed as
a generalization of the familiar fact that
the Lie bracket of symplectic vector fields is
a Hamiltonian vector field.

\section{FWX resolution of the algebra of local Hoch\-schild
co\-chains}
In this section we recall the construction of
the Fedosov-Weinstein-Xu
(FWX) resolution \cite{WX} of the algebra of local Hochschild cochains
of the quantum ring of functions on a symplectic manifold.
A classical analogue of this resolution was used in papers
\cite{CEFT} and \cite{FTC} to prove the formality theorems
for Hochschild (co)chains of the algebra of functions
on an arbitrary smooth manifold. We start with
a brief reminder of the Fedosov
deformation quantization.

\subsection{Reminder of the Fedosov deformation quantization}
As above, $M$ stands for an even dimensional smooth manifold
endowed with a symplectic form
$\om=\om_{ij}(x)dx^i\wedge dx^j$ and
$\om^{ij}(x)\pa_{x^i}\wedge \pa_{x^j}$ denotes the respective
Poisson tensor
$$
\om^{ik}(x)\om_{kj}(x)=\de^i_j\,.
$$

Following Fedosov \cite{Fedosov} we introduce
the Weyl algebra bundle over the manifold $M$.
Sections of this bundle are the following formal sums

\begin{equation}
\label{sect}
a=a(x,\h,y)=\sum_{p\ge 0,\, k\in \bbZ} \h^k a_{k;i_1 \dots i_p}(x)y^{i_1} \dots y^{i_p}\,,
\end{equation}
where the summation in $k$ is bounded below, $a_{k; i_1 \dots i_p}(x)$
are symmetric covariant $\bbC$-valued tensors, and $y^i$ are
fiber coordinates on the tangent bundle $TM$\,.
An associative multiplication of the sections (\ref{sect})
is defined with the help of the Poisson tensor $\om^{ij}$
as follows

\begin{equation}
\label{circ}
a\circ b(x,\h,y)=exp\left( \frac{\h}2\om ^{ij}(x)\frac \partial {\partial
y^i}\frac \partial {\partial z^j}\right) a(x,\h, y)\,
b(x,\h,z)|_{z=y}\,.
\end{equation}

The bundle $\W$ viewed as a sheaf of algebras over $\bbC$
is naturally filtered with respect
to the degree of monomials $2[\h]+[y]$ where
$[\h]$ is a degree in $\h$ and $[y]$ is a degree in $y$

\begin{equation}
\label{filtr}
\begin{array}{c}
\displaystyle
\dots \subset \W^1  \subset \W^0  \subset \W^{-1}
\dots \subset \W\,, \\[0.3cm]
\displaystyle
\G(X,\W^m) = \{ a = \sum_{2k+p\ge m}
\h^k a_{k;i_1 \dots i_p}(x)y^{i_1} \dots y^{i_p}\}
\end{array}
\end{equation}
This filtration defines the $2[\h]+[y]$-adic topology in the
algebra of section $\G(X,\W)$ over any open subset $X\subset M$\,.

Let us remark that the space $\Om^{\bul}(M, \W)$ of smooth exterior
forms with values in $\W$ is naturally a graded associative algebra
with the product induced by (\ref{circ}) and the
following graded commutator
$$
[a,b]= a \circ b - (-)^{q_a q_b} b\circ a\,,
$$
where $q_a$ and $q_b$ are exterior degrees of
$a$ and $b$\,, respectively.
The filtration of $\W$ (\ref{filtr}) gives us a filtration
of the algebra $\Om^{\bul}(M, \W)$

\begin{equation}
\label{filtrOm}
\dots \subset \Om(M,\W^1) \subset
\Om(M,\W^0) \subset \Om(M,\W^{-1})
\dots \subset \Om(M,\W)\,,
\end{equation}
which similarly defines
the $2[\h]+[y]$-adic topology in the algebra
$\Om(M, \W)$ of exterior forms\,.

Following Fedosov \cite{Fedosov} we introduce the following
derivation of the algebra $\Om(M,\W)$
\begin{equation}
\label{del}
\de= dx^i \frac{\pa}{\pa y^i} \,:\, \Om^{\bul}(M,\W) \mapsto
\Om^{\bul+1}(M,\W)\,,
\qquad \de^2=0\,.
\end{equation}

The cohomology of the differential $\de$ are described by
the following
\begin{pred}
\label{coh-de}
$$
H^{0}(\OmW,\de) = C^{\infty}(M)((\h))\,, \qquad
H^{\ge 1} (\OmW, \de)=0\,.
$$
\end{pred}
{\bf Proof.} It is not hard to guess the map
$$
\si\,:\, \OmW \mapsto C^{\infty}(M)((\h))
$$
from the complex $(\OmW, \de)$ onto its subcomplex
$(C^{\infty}(M)((\h)),0)$ and a contracting homotopy $\de^{-1}$
$$
\de^{-1} \,:\,\Om^{\bul}(M,\W)\mapsto \Om^{\bul-1}(M,\W) \,,
$$
between the map $\si$ and the identity map. Namely,

\begin{equation}
\label{sigma}
\si (a)= a \Big |_{y=0,~dx=0}\,, \qquad
a\in \Om^{\bul}(M,\W)\,.
\end{equation}
and

\begin{equation}
\delta^{-1}a=y^k i\left( \frac \partial {\partial x^k}\right)
\int\limits_0^1 a(x,\h, ty,tdx)\frac{dt}t,  \label{del-1}
\end{equation}
where
$i(\partial /\partial x^k)$ denotes the contraction
of an exterior form with the vector field
$ \partial /\partial x^k$\,, and $\delta ^{-1}$ is extended to
$\G(\W)$ by zero.

The desired contracting property
\begin{equation}
a=\sigma(a) +\delta \delta ^{-1}a + \delta ^{-1}\delta a\,,
\qquad \forall~a\in\OmW
 \label{Hodge}
\end{equation}
can be checked by a straightforward
computation. $\Box$\\
We remark that the operator $\de^{-1}$ turns
out to be nilpotent
\begin{equation}
\label{de-1-nilp}
(\de^{-1})^2=0\,.
\end{equation}

Recall that on any symplectic manifold there exists an
affine torsion-free connection $\n_i$ which is compatible with
the symplectic structure \cite{GRSh}
$$
\n_i \om_{jk}(x)=0\,.
$$
This connection allows us to define the
following derivation of the algebra $\OmW$
\begin{equation}
\label{nab}
\n = dx^i\frac{\pa}{\pa x^i} -
dx^i\G^j_{ik}(x) y^k\frac{\pa}{\pa y^j}
\,:\,\Om^{\bul}(M,\W)\mapsto \Om^{\bul+1}(M,\W)\,,
\end{equation}
where $\G^j_{ik}(x)$ are the respective Christoffel symbols.
The compatibility of the connection with the symplectic structure
implies that (\ref{nab}) is a derivation of
the product (\ref{circ})\,.

Since the connection $\n_i$ is torsion-free
the derivation $\n$ anti-commutes with $\de$
\begin{equation}
\label{n-de=de-n}
\n \de + \de \n\,.
\end{equation}

In general the derivation $\n$ is not nilpotent as $\de$. Instead we
have the following expression for $\n^2$

\begin{equation}
\label{nab-sq}
\n^2 a = \frac1{\h}[R , a]
\,:\, \Om^{\bul}(M,\W) \mapsto \Om^{\bul+2}(M,\W)\,,
\end{equation}
where
$$
R= -\frac14 dx^i dx^j \om_{km}(R_{ij})^m_l(x) y^k y^l \,,
$$
and $(R_{ij})^k_l(x)$ is the standard Riemann curvature tensor
of the connection $\n_i$.

Let us consider the following derivation of the
algebra $\OmW$ (connection on the bundle $\W$)
\begin{equation}
\label{DDD}
D=\n -\de + \frac{1}{\h}[\,r,\,\bul\,]\,,
\end{equation}
where $r$ is a smooth $1$-form with values in
$\W^{3}$. Using equations (\ref{n-de=de-n}) and (\ref{nab-sq})
we get that for any $a\in \OmW$
\begin{equation}
\label{D-sq}
D^2 a = \frac1{\h}[R-\de r + \n r + \frac1{\h}r \circ r , a]\,.
\end{equation}

Since the center of the Weyl algebra is
$\bbC((\h))\subset \bbC[[y^1, \dots , y^{2n}]]((\h))$
the derivation $D$ is nilpotent (the connection
(\ref{DDD}) is flat) if and only if
\begin{equation}
\label{iff}
R-\de r + \n r + \frac1{\h}r \circ r \in \Om^2(M)((\h))\,.
\end{equation}
Furthermore, since $r\in \Om^1(M,\W^3)$ we have
\begin{equation}
\label{iff1}
R-\de r + \n r + \frac1{\h}r \circ r \in \h\Om^2(M)[[\h]]\,.
\end{equation}
Let us say that

\begin{equation}
\label{c-class}
R-\de r + \n r + \frac1{\h}r \circ r = \sum_{k\ge 1}\h^{k}\om_k\,.
\end{equation}
Using the properties of $\n$ and $\de$ and the Bianchi
identities for the Riemann curvature tensor $(R_{ij})^k_l$
$$
\de R=0\,, \qquad  \n R =0\,.
$$
we derive that
$$
D (R-\de r + \n r + \frac1{\h}r \circ r ) = 0\,.
$$
The latter is equivalent to the fact that all the forms
$\om_k$ in the right hand side of (\ref{c-class}) are
closed with respect to the De Rham differential $d$\,.

It turns out that for any series of closed forms

\begin{equation}
\label{class}
\Om = \sum_{k\ge 1}\h^{k}\om_k
\end{equation}
one can construct an element $r\in \Om^1(M, \W^3)$
satisfying (\ref{c-class}). Namely, theorem $5.3.3$
in \cite{Fedosov1} implies that

\begin{teo}[Fedosov]
For any formal series (\ref{class}) of closed forms on $M$
one can construct an element
\begin{equation}
\label{rrr}
r= \sum_{k\ge 0,\, p\ge 1} dx^l  \h^k r_{k; l, i_1 \dots i_p}(x)y^{i_1} \dots y^{i_p}
\in \Om^{1}(M, \W^{3})
\end{equation}
satisfying (\ref{c-class}) and the
condition $\de^{-1}r=0$\,. The operation
\begin{equation}
\label{Fed-diff}
D=\n -\de + \frac{1}{\h}[\,r,\,\bul\,]
\end{equation}
associated to the element (\ref{rrr}) is
a nilpotent derivation of
the algebra $\OmW$\,.
\end{teo}
{\bf Remark.} The element $r$ in the above theorem can
be obtained by iterating the following equation
\begin{equation}  \label{iter}
r=\delta ^{-1}(R-\Om) + \de^{-1}(\nabla r+
\frac 1{\h }r\circ r)\,.
\end{equation}
Since the operator $\nabla $ preserves filtration (\ref{filtrOm})
and $\delta ^{-1}$ raises it by $1$, the iteration of (\ref{iter})
converges in the topology (\ref{filtrOm}) and defines the unique element
$r \in \Om^{1}(M, \W^{3})$\,.

Furthermore, the respective generalization of theorem $5.2.4$
in \cite{Fedosov1} says that

\begin{teo}[Fedosov]
\label{ONA}
Iterating the equation

\begin{equation}
\label{iter-a}
\tau (a)= a + \de^{-1}(\n\tau (a) + \frac1{\h}[r,\tau(a)])
\end{equation}
for $a\in C^{\infty}(M)((\h))$ one constructs an
isomorphism
$$
\tau\,:\,C^{\infty}(M)((\h)) \mapsto ker\,D\cap \G(M,\W)
$$
from the vector space $C^{\infty}(M)((\h))$ to the
vector space of horizontal sections $Z^0(\OmW, D)$ of the
connection (\ref{Fed-diff}). For any $a\in C^{\infty}(M)((\h))$
$$
\si (\tau (a)) = a
$$
and the multiplication on $C^{\infty}(M)((\h))$ induced by
the $\circ$-product (\ref{circ}) via the isomorphism $\tau$

\begin{equation}
\label{star}
a \ast  b= \si(\tau(a) \circ \tau(b)), \qquad
a, b\in C^{\infty}(M)((\h))
\end{equation}
is a star-product associated to
the Poisson bracket $\om^{ij}$\,. ~$\Box$
\end{teo}
Notice that since the operator $\nabla $ preserves the filtration (\ref{filtrOm})
and $\delta ^{-1}$ raises it by $1$, the iteration of (\ref{iter-a}) converges in the
topology (\ref{filtrOm}) and define the unique element
$\tau(a)\in Z^0(\OmW,D)$ for any $a\in C^{\infty}(M)((\h))$.

One can easily observe that the series (\ref{class}) of closed
forms enters the construction of the star-product (\ref{star}).
In fact it is not hard to show that the equivalence class of the star-product
(\ref{star}) depends only on the cohomology class of
the series (\ref{class}). In other words, if $\ast$ and
$\ast'$ are Fedosov star-products corresponding to
series $\Om$ and $\Om'$ representing the same cohomology
class in $H^2(M,\bbC)((\h))$ then $\ast$ is equivalent to
$\ast'$ in the sense of (\ref{P-equiv}).
Furthermore, if two series $\Om$ and $\Om'$
define distinct cohomology classes in $H^2(M,\bbC)((\h))$
then the corresponding Fedosov star-products $\ast$ and
$\ast'$ are not equivalent.

The following result of P. Xu \cite{Xu} shows that
the above construction allows us to get a star-product
from any equivalence class of the star-products on $M$
\begin{teo}[P. Xu, \cite{Xu}]
\label{Xu}
Any star-product on a symplectic (smooth real) manifold
$M$ is equivalent to some Fedosov star-product.~$\Box$
\end{teo}
If a star-product $*$ on $M$ is equivalent to the Fedosov
star-product corresponding to the series (\ref{class})
the cohomology class of the
combination
\begin{equation}
\label{repres}
\Om_F=\frac1{\h}(-\om + \Om)\in \frac1{\h}Z^2_d(\Om(M)[[\h]])
\end{equation}
is called the Fedosov class of a star-product $*$.
We refer to $\Om_F$ entering the construction of
the star-product as {\it the Fedosov representative}.

\subsection{Double complex of fiberwise Hoch\-schild co\-chains.}
In this section we denote by $\bbA$ the algebra
of functions $C^{\infty}(M)((\h))$ with the
star-product $\ast$ (\ref{star}) associated
to the Fedosov data $(M,\n,\Om_F)$.

We now turn to the definition of formal fiberwise
Hochschild cochains on $\G(\W)$

\begin{defi}
A bundle $\cC^k $ of formal fiberwise
Hochschild cochains of degree $k$ is a bundle over $M$ whose sections
are $C^{\infty}(M)$-polylinear
maps $\mP : \bigotimes^{k} \G(\W) \mapsto \G(\W)$
continuous in the adic topology (\ref{filtr}).
\end{defi}
Any such map $\mP\in\G(\cC^k)$ can be uniquely
represented in the form of the following formal series
\begin{equation}
\label{operr}
\mP =\sum_{\al_1 \dots \al_k}\sum_{m,p=0}^{\infty}
\h^m \mP^{\al_1\dots \al_k}_{m, i_1\dots i_p}(x)y^{i_1}
\dots y^{i_p} \frac{\pa}{\pa y^{\al_1}}\otimes  \dots \otimes \frac{\pa}{\pa
y^{\al_k}}\,,
\end{equation}
where the summation in $m$ is bounded below,
$\al$'s are multi-indices $\al={j_1\dots j_l}$\,,
$$
\frac{\pa}{\pa y^{\al}}=\frac{\pa}{\pa y^{j_1}} \dots \frac{\pa}{\pa
y^{j_l}}\,,
$$
and the tensors $\mP^{\al_1\dots \al_k}_{m; i_1\dots i_p}(x)$ are
symmetric in covariant indices $i_1,\dots, i_p$\,.

Extending the above definition by allowing the cochains
to be inhomogeneous we define the total bundle $\cC$ of
formal fiberwise Hochschild cochains as a direct sum

\begin{equation}
\label{cal-D}
\cC =\bigoplus_{k=0}^{\infty} \cC^k\,, \qquad
\cC^{0}=\W\,.
\end{equation}

The space $\OmC$ of smooth exterior forms with values in
the bundle $\cC$ acquires a natural associative product
induced by the fiberwise cup-product (\ref{cup-W}) in
the space of Hochschild cochains of the formal
Weyl algebra. Here we also call this product a
cup-product. So, given $\mP_1\in \Om(M, \cC^{k_1})$ and
$\mP_2\in \Om(M, \cC^{k_2})$ their
cup-product is defined by
\begin{equation}
\label{cup}
\mP_1 \cup \mP_2(a_1, a_2,\dots, a_{k_1+k_2}) =
\mP_1 (a_1, \dots, a_{k_1}) \circ
\mP_2 (a_{k_1+1},\dots, a_{k_1+k_2})\,,
\end{equation}
where $a_i$ are arbitrary smooth sections of
the Weyl algebra bundle $\W$\,.

The space $\OmC$ also acquires a graded Lie algebra
structure induced by the fiberwise
Gerstenhaber bracket (\ref{Gerst-W}).
For two homogeneous elements $\mP_1\in \Om(M,\cC^{k_1+1})$ and
$\mP_2\in \Om(M, \cC^{k_2+1})$ the bracket
is defined as follows

$$[\mP_1, \mP_2]_G(a_0,\,\dots, a_{k_1+k_2})=$$
\begin{equation}
\label{Gerst}
\sum_{i=0}^{k_1}(-)^{ik_2}
\mP_1(a_0,\,\dots , \mP_2 (a_i,\,\dots,a_{i+k_2}),\, \dots, a_{k_1+k_2})
\end{equation}
$$
-(-)^{k_1k_2} (1 \leftrightarrow 2)\,, \qquad a_j \in \G(M,\W)\,.
$$

The fiberwise $\circ$ product in $\G(\W)$ gives
us the fiberwise Hochschild
differential (\ref{Hoch-diff})

\begin{equation}
\label{Hoch-d}
(\pa \mP) (a_0, a_1, \dots, a_k) = (-)^q (a_0 \circ \mP(a_1, \dots, a_k)
- \mP(a_0\circ a_1, a_2, \dots, a_k)+\dots +
\end{equation}
$$
(-1)^{k}\mP(a_0, \dots, a_{k-2}, a_{k-1}\circ a_k)+
(-1)^{k+1} \mP(a_0, a_1, \dots, a_{k-1})\circ a_k)\,,
$$
$$
\pa\,:\, \Om^{q}(M,\cC^k) \mapsto \Om^{q}(M, \cC^{k+1})\,,
$$
which is a derivation of both the cup-product (\ref{cup}) and
the Gerstenhaber bracket (\ref{Gerst}) by (\ref{cup-der})
(\ref{G-der-pa}).

The Fedosov differential (\ref{DDD}) can be naturally extended
to the vector space $\OmC$ via the formula
$$
(D \mP) (a_1, \dots, a_k) =
$$
\begin{equation}
\label{DDD-coch}
 D \mP (a_1, \dots, a_k) -(-)^{q}
(\mP (Da_1, a_2 \dots, a_k) + \dots + \mP (a_1, a_2 \dots, D
a_k))\,,
\end{equation}
where $\mP\in \Om^q(M,\cC^k)$ and $a_i$ are arbitrary
smooth sections of the Weyl algebra bundle $\W$\,.
It is not hard to see that the differential (\ref{DDD-coch})
is a derivation of the cup-product (\ref{cup}) and
the Gerstenhaber bracket (\ref{Gerst})\,. The differential
(\ref{DDD-coch}) also anti-commutes with
the fiberwise Hochschild differential $\pa$
$$
D \, \pa + \pa \, D = 0
$$
since $D$ is a derivation of the $\circ$-product
(\ref{circ}).

Thus we arrive at the double complex
$(\Om^{\bul}(M, \cC^{\bul}), D, \pa)$, the total space
of which is a differential graded associative algebra (DGAA)
and also a differential graded Lie algebra (DGLA).

Due to propositions \ref{homot-Hoch} and \ref{equiv}
given in Appendix B the fiberwise Hochschild
differential $\pa$ has vanishing higher cohomology
\begin{equation}
\label{ge1=0}
H^{\ge 1}(\Om(M, \cC^{\bul}),\pa) = 0\,.
\end{equation}
Therefore for any $q\ge 0$ the natural embedding

\begin{equation}
\label{vr-q}
\mE^q\,:\,\Om^q(M)((\h)) \hookrightarrow \Om^q(M, \cC^{0})
\end{equation}
induces a quasi-isomorphism of complexes
$(\Om^q(M)((\h)),0)$ and $(\Om^q(M,\cC^{\bul}), \pa)$.
Furthermore, since the Fedosov differential restricted
to $\Om(M)((\h))$ coincides with the De Rham differential $d$
the inclusion $\mE\,:\,\Om(M) \hookrightarrow \Om(M, \cC^{0})$
also induces a morphism from the De Rham complex
$(\Om(M)((\h)), d)$ to the total complex
$(\Om(M, \cC^{\bul}), D+\pa)$.
One can easily see that $\mE$ is
compatible with the cup-products and therefore
$\mE$ is a morphism of DG associative
algebras

\begin{equation}
\label{vrooo}
\mE\,:\,(\Om(M)((\h)), d, \wedge)
\hookrightarrow (\Om(M, \cC), D+\pa, \cup)\,.
\end{equation}
Due to (\ref{ge1=0}) and (\ref{vr-q}) we have that
\begin{pred}
\label{AGA}
The spectral sequence of the
double complex $(\Om^{\bul}(M, \cC^{\bul}), D, \pa)$\,,
associated to the filtration by the degree of
fiberwise Hochschild cochains degenerates at $E_1$ and
the map $\mE$ (\ref{vrooo}) is a quasi-isomorphism  of DG associative
algebras. $\Box$
\end{pred}

We will prove the desired statements of theorems
$1$ and $2$ by computing the total cohomology
of double complex $(\Om^{\bul}(M, \cC^{\bul}), D, \pa)$
using the spectral sequence
associated to the filtration by the exterior degree
in $\Om^{\bul}(M, \cC^{\bul})$\,. But before doing
this, we need to extend the operations
introduced on the vector space $\OmW$ in the previous subsection
to the vector space $\Om(M, \cC)$\,. A proof of the following proposition
is straightforward
\newpage
\begin{pred}~\\[-0.5cm]
\label{EEE}
\begin{enumerate}
\item The nilpotent derivation $\de$ (\ref{del}) and
the covariant derivative $\n$ (\ref{nab}) are extended
to $\Om(M,\cC)$ in the following natural manner

\begin{equation}
\label{del-OmC}
\de \mP = [dx^i\frac{\pa}{\pa y^i}, \mP]_G\,,
\qquad
\n \mP = dx^i \frac{\pa}{\pa x^i} \mP -
[\,dx^i \G^j_{ik} y^k\frac{\pa}{\pa y^j}\,,\,\mP]_G\,,
\end{equation}
where
$\displaystyle dx^i \G^j_{ik} y^k\frac{\pa}{\pa y^j}$ is
viewed locally as an element of $\Om^{1}(\,\bul\,, \cC^1)$\,.
With this definition $\n$ is globally defined, $\de$
is nilpotent and equation (\ref{n-de=de-n}) holds.

\item The component-wise extension of
the map $\si$ (\ref{sigma})

\begin{equation}
\label{sigma-OmC}
\si \mP = \mP\Big|_{y^i=dx^i=0}\,:\,
\OmC \mapsto Z^0_{\de}(\OmC)
\end{equation}
defines a projection onto the
kernel $Z^0_{\de}(\OmC)= ker\,\de \cap \G(\cC)$
of $\de$ in $\G(\cC)$\,.
With this definition of $\si$ and the component-wise
definition of $\de^{-1}$ (\ref{del-1}) equations
(\ref{Hodge}) and (\ref{de-1-nilp}) hold in
$\Om(M,\cC)$\,.

\item The Fedosov differential on $\OmC$ (\ref{DDD-coch})
can be rewritten in the form

\begin{equation}
\label{DDD-cC}
D \mP = \n \mP - \de \mP + \frac1{\h}
[\pa r ,\mP ]_G\,, \qquad \forall~\mP\in \OmC\,,
\end{equation}
where $r$ (\ref{rrr}) is viewed as an element in
$\Om^1(M,\cC^0)$\,. $\Box$
\end{enumerate}
\end{pred}

The following proposition shows that the complex
$\Om^{\bul}(M, \cC^{\bul})$ is a resolution of
the complex $C^{\bul}_{loc}(\bbA)$ of local Hochschild
cochains of $\bbA$
\begin{pred}
\label{main}
Fedosov differential (\ref{DDD-coch})
has vanishing higher cohomology
\begin{equation}
\label{higher=0}
H^{\ge 1} (\Om(M, \cC^{\bul}),D)=0
\end{equation}
and
\begin{equation}
\label{0=awesome}
H^0 (\Om(M, \cC^{\bul}),D)=
C_{loc}^{\bul}(\bbA)
\end{equation}
as DG associative algebras and as DG Lie algebras.
\end{pred}
{\bf Proof.}
To prove (\ref{higher=0}) we pick up
$\mP\in \Om^{\ge 1}(M, \cC)$ which is closed
with respect to the Fedosov differential  $D$
$$
D \mP=0
$$
and observe that the recurrent procedure

\begin{equation}
\label{iter-mQ}
\mQ = -\de^{-1} \mP + \de^{-1}
(\n \mQ + \frac{1}{\h} [\pa r, \mQ]_G)
\end{equation}
converges in $2[\h]+[y]$-adic topology to
an element $\mQ\in \Om(M, \cC)$ such that
\begin{equation}
\label{mQ-aux1}
\mQ \Big|_{y=0}=0\,.
\end{equation}
Due to equations (\ref{Hodge}), (\ref{de-1-nilp}),
extended to $\OmC$ by proposition \ref{EEE},
and equation (\ref{mQ-aux1})
$$
\de^{-1}\mQ=0\,,
$$
and
\begin{equation}
\label{mQ-aux}
\de^{-1}(D\mQ -\mP)=0\,.
\end{equation}
Let us denote $D\mQ-\mP$ by $\mK$.
Using (\ref{Hodge}) once again we get
that $\mK$ satisfies
\begin{equation}
\label{mK-prop}
\mK= \de^{-1} (\n\mK +\frac1{\h}[\pa r,\mK]_G)
\end{equation}
since $\mK\in \Om^{\ge 1}(M,\cC)$ and
$D \mK =0$.

Equation (\ref{mK-prop}) has the unique vanishing solution
since $\de^{-1}$ raises the degree in $y$\,. Hence (\ref{higher=0})
is proven.

To prove the second assertion we mention
the following property of the map $\tau$
(\ref{iter-a})
\begin{equation}
\label{der-tau-a}
\frac{\pa}{\pa y^{i_1}}
\dots \frac{\pa}{\pa y^{i_k}} \tau (a)\Big|_{y^i=0}(x,\h) =
\pa_{x^{i_1}} \dots \pa_{x^{i_k}} a(x,\h) + ~{\rm
lower~order~derivatives~in}~x\,,
\end{equation}
$$\forall~ a\in C^{\infty}(M)((\h))\,.$$
Thus using the map $\tau$ we can identify
the vector space $C^{\bul}_{loc}(\bbA)$
of local Hochschild cochains of $\bbA$ and
the vector space $Z^0_{\de}(\OmC)= ker\, \de \cap
\G(\cC)$\,. An isomorphism from the latter
space to the former one is given by the formula

\begin{equation}
\label{mu}
(\mu \mP) (a_1, \dots, a_k) =
\si\mP (\tau (a_1), \dots, \tau(a_k))\,,
\end{equation}
where $\mP\in Z^0_{\de}(\OmC)$ and
$a_i$ are elements in $C^{\infty}(M)((\h))$\,.
It follows form (\ref{mu}) that the star-product
(\ref{star}) is the image of $\mu$, namely
\begin{equation}
\label{mu-circ}
\ast = \mu (\circ)\,,
\end{equation}
where $\circ$ is viewed as an element in
$\G(\cC^2)$ and $\ast\in C^2_{loc}(\bbA)$\,.

An appropriate modification of theorem \ref{ONA}
(or theorem $3$ in \cite{CEFT}) enables us to conclude
that iterating the equation
\begin{equation}
\label{iter-mP}
\al(\mP)= \mP + \de^{-1}(\n \al(\mP) +
\frac1{\h}[\pa r,\al(\mP)]_G)
\end{equation}
for any $\mP\in Z^0_{\de}(\OmC)$ one constructs an
isomorphism
$$
\al\,:\,Z^0_{\de}(\OmC)\, \widetilde{\longrightarrow}\,
Z^0_D(\OmC)\,,
$$
where $Z^0_D(\OmC)= ker\, D\cap \G(\cC)$\,. It is obvious
that
$$
\si (\al (\mP))= \mP \qquad \forall~\mP \in Z^0_{\de}(\OmC)\,.
$$
Notice that the element $\circ\in \G(\cC^2)$
remains unchanged under $\al$
$$
\al(\circ) = \circ\,.
$$

The second claim of the proposition will follow
if we prove that the map
\begin{equation}
\label{beta}
\beta(\mP)= \mu (\si(\mP)) \,:\,
Z^0_D(\Om(M, \cC^{\bul}))\mapsto C^{\bul}_{loc}(\bbA)
\end{equation}
is an isomorphism of DG associative algebras and
DG Lie algebras.
We already know that $\beta$ is an one-to-one and
onto. Thus we have to prove the compatibility
with the algebraic operations $\cup$, $[,]_G$, and
$\pa$. To do this we observe that for any $k\ge 0$
the map $\beta$ is given by the
formula\footnote{For $k=0$ the map $\beta$
just coincides with $\si$}
\begin{equation}
\label{beta1}
(\beta \mP)(a_1, \dots, a_k)=
\si(\mP (\tau(a_1), \dots, \tau (a_k)))\,.
\end{equation}
The compatibility with the cup-product (\ref{cup})
follows from the line of equations
$$
(\beta(\mP_1 \cup \mP_2)) (a_1, \dots, a_{k_1+k_2})=
\si (\mP_1(\tau(a_1), \dots, \tau(a_{k_1}))\circ
\mP_2(\tau(a_{k_1+1}), \dots, \tau(a_{k_1+k_2})))=
$$
$$
\si (\mP_1(\tau(a_1), \dots, \tau(a_{k_1})))\ast
\si (\mP_2(\tau(a_{k_1+1}), \dots, \tau(a_{k_1+k_2})))=
(\beta \mP_1)\cup (\beta \mP_2) (a_1, \dots, a_{k_1+k_2})\,,
$$
where $a_1, \dots, a_{k_1+k_2}$ are arbitrary elements
in $\bbA$, $\mP_1\in Z^0_D(\Om(M, \cC^{k_1}))$, and
$\mP_2 \in Z^0_D(\Om(M, \cC^{k_2}))$\,.

To prove the compatibility with the Gerstenhaber bracket
we observe that for any $a_1, \dots, a_k\,\in \bbA$
and $\mP\in Z^0_D(\Om(M, \cC^k))$
\begin{equation}
\label{OGO}
\tau ((\beta \mP)(a_1, \dots, a_k))=
\mP(\tau(a_1), \dots, \tau (a_k))\,.
\end{equation}
The latter is proven as follows. Both right and
left hand sides of (\ref{OGO}) are $D$-closed
elements of $\G(\W)$. By (\ref{beta1}) $\si$ of
the left hand side equals to $\si$ of the right hand
side. Hence (\ref{OGO}) follows from theorem \ref{ONA}\,.

Using (\ref{OGO}) we derive that for any
$\mP_1\in Z^0_D(\Om(M, \cC^{k_1}))$,
$\mP_2\in Z^0_D(\Om(M, \cC^{k_2}))$, and
$a_1, \dots, a_{k_1+k_2-1}\,\in \bbA$
$$
( \beta \mP_1(\dots, \mP_2(\dots), \dots ))
(a_1, \dots, a_{k_1+k_2-1})=
$$
$$
\si \mP_1(\tau(a_1), \dots, \mP_2(\tau(a_i), \dots,
\tau(a_{i+k_2-1})), \dots, \tau(a_{k_1+k_2-1}) )=
$$
$$
\si \mP_1(\tau(a_1), \dots, \tau ((\beta\mP_2)(a_i, \dots,
a_{i+k_2-1})), \dots, \tau(a_{k_1+k_2-1}) )=
$$
$$
(\beta \mP_1)(a_1, \dots, (\beta \mP_2)(a_i, \dots,
a_{i+k_2-1}), \dots, a_{k_1+k_2-1})\,,
$$
where $\mP_1(\dots, \mP_2(\dots), \dots )$ denotes
the substitution of $\mP_2$ on the place of the $i$-th
argument of $\mP_1$\,. Thus the compatibility of (\ref{beta})
with the Gerstenhaber bracket is proven.

The compatibility of (\ref{beta}) with the Hochschild differential
follows from the compatibility of (\ref{beta})
with the Gerstenhaber bracket, and
equations (\ref{mu-circ}), (\ref{pa}).
Thus the proposition is proven. $\Box$

~\\
{\bf Proofs of theorems $1$ and $2$.}
It turns out that most of work is already
done. Due to theorem \ref{Xu} of P. Xu \cite{Xu}
we can safely assume that $\ast$ is the
star-product associated to the Fedosov data
$(M,\n,\Om_F)$\,. Then it is clear that the
first claim of theorem $1$ would follow from
theorem $2$\,.

Proposition \ref{main} implies that the
spectral sequence of the double complex
$(\Om^{\bul}(M,\cC^{\bul}),D,\pa)$ associated
to the filtration by the exterior degree
degenerates at $E_1$ and the total cohomology
of $(\Om^{\bul}(M,\cC^{\bul}),D,\pa)$
$$
H^{\bul}(\Om^{\bul}(M,\cC^{\bul}), D+\pa) =
H^{\bul}(Z^0_D(\Om^{\bul}(M,\cC^{\bul})), \pa)
$$
both as graded associative algebras and
as graded Lie algebras. Combining the statements of
propositions \ref{AGA} and \ref{main} we get
that $(\Om^{\bul}(M,\cC^{\bul}), D+\pa, \cup)$
is the desired DG associative algebra
$K^{\bul}$ and the quasi-isomorphisms
in question
\begin{equation}
\label{leftright}
(\Om^{\bul}(M)((\h)), d, \wedge)
\,\stackrel{\mE}{\longrightarrow}\,
(\Om^{\bul}(M,\cC^{\bul}), D+\pa, \cup)
\, \stackrel{\mD}{\longleftarrow}\,
(C^{\bul}_{loc}(\bbA), \pa, \cup)
\end{equation}
are $\mE$ (\ref{vrooo}) and $\mD= \beta^{-1}$
(\ref{beta}).
The naturality of $K^{\bul}$, $\mD$ and $\mE$ with respect
to the Fedosov data is evident.

While the equivariance property (\ref{al-equiv}) of $\mE$
is obvious from the construction, the equivariance property
(\ref{al-equiv}) of $\mD$ follows essentially
from the fact that any
diffeomorphism acts on the fiber variables $y^i$ of the
tangent bundle $TM$ by linear transformations.
Thus theorem $2$ is proven.

To prove the statement about the Gerstenhaber bracket
in theorem $1$ we observe that due to propositions \ref{homot-Hoch} and
\ref{equiv} in Appendix B any cocycle in
$(\Om^{\bul}(M,\cC^{\bul}), D+\pa)$ is cohomologically
equivalent to a cocycle in $\G(\cC^{0}))$\,.
But the restriction of the Gerstenhaber bracket to
$\cC^0$ is vanishing by definition (\ref{Gerst-W}).
This completes the proofs of both theorems. $\Box$

\section{Concluding remarks}
In this section we discuss applications and variations
of theorems $1$ and $2$\,.

First, using theorems $1$ and $2$\,, one can
easily prove an equivariant version of
Xu's theorem \cite{Xu}

\begin{cor}
If $G$ is a group acting on $M$ by symplectomorphisms
and $M$ admits a $G$-invariant connection $\n$ then
any $G$-invariant star-product is equivalent
to some $G$-invariant Fedosov star-product.
If $G$ is finite or compact then
the equivalence can be established by a $G$-invariant
operator. $\Box$
\end{cor}
Second, in some cases it is instructive to know
the Hochschild cohomology of the algebra
$\bbA_0= C^{\infty}(M)[[\h]]$ of formal Taylor power
series with multiplication $\ast$. By keeping track of
negative degrees in $\h$ in our construction
one can easily prove that
\begin{pred}
The graded associative algebras
$$
(HH^{\bul}_{loc}(\bbA_0), \cup)
$$
and
$$
H^{\bul}(\Om(M)[[\h]], \h\,d, \wedge)
$$
are isomorphic. $\Box$
\end{pred}

Third, natural algebraic geometric versions
of theorems $1$ and $2$ hold for a smooth affine
algebraic variety $X$ (over $\bbC$).
These versions immediately follow from
Grothendieck's theorem on De Rham cohomology
of an affine variety and the fact that any smooth affine algebraic variety
admits an algebraic connection on the holomorphic tangent
bundle $T_{hol}X$\,.
\begin{pred}
If $X$ is a smooth affine algebraic variety over $\bbC$
endowed with an algebraic symplectic form $\om$,
and $\bbA =(\cO(X)((\h)), \ast)$ is the corresponding
quantum ring of functions then the graded
associative algebra
$$
(HH^{\bul}(\bbA), \cup)
$$
of Hochschild cohomology of $\bbA$
is isomorphic to the graded associative
algebra
$$
(H^{\bul}_{DR}(X), \cup)
$$
of De Rham cohomology of $X$.
The Gerstenhaber bracket on $HH^{\bul}(\bbA)$ is
vanishing. $\Box$
\end{pred}
It is worth mentioning that in this algebraic
geometric case the complex of local Hochschild
cochains is quasi-isomorphic to the complex
of all Hochschild cochains. For this reason,
the analogue of theorem $1$ is formulated
for the genuine Hochschild cohomology
of the quantum ring $\bbA$\,.

Notice that using the latter proposition only for
the Fedosov star-products and rearranging the arguments,
one can prove an algebraic geometric version
of Xu's theorem \cite{Xu} (theorem \ref{Xu}) for
any smooth symplectic affine algebraic variety
over $\bbC$\,.

Finally, we would like to mention paper \cite{FFS} in which
the authors propose an explicit expression for the canonical
trace in deformation quantization of a symplectic manifold.
We suspect that their formula for the trace can be
obtained with the help of quasi-isomorphisms between
the bar resolution and the Koszul resolution of
Weyl algebra (see proposition \ref{K-bar-K})
and the origin of the integrals
over the configuration space of ordered points on a circle
in their formula is the result of multiple applications
of a contracting operator similar to (\ref{homot-H}).

{\bf Acknowledgment.}
This paper arose from questions of my advisor
Pavel Etingof to whom I express my sincere thanks.
I also acknowledge Pavel for his valuable
criticisms concerning the first version of the manuscript.
I would like to thank Alexander Braverman, Simone Gutt,
Lars Hesselholt, Richard Melrose, Alexei Oblomkov, and
Dmitry Tamarkin for useful discussions.
I am grateful to D. Silinskaia for criticisms
concerning my English.
The work is partially supported by the
NSF grant DMS-9988796, the Grant
for Support of Scientific Schools NSh-1999.2003.2,
the grant INTAS 00-561 and the grant CRDF
RM1-2545-MO-03.

\section{Appendix A. Algebraic structures on Hoch\-schild
co\-chains.}
\renewcommand{\theequation}{A.\arabic{equation}}
Let $A$ be an associative unital algebra over a field of characteristic
zero. By definition, Hochschild cohomology \cite{Loday} of $A$ is
\begin{equation}
\label{Hoch}
HH^{\bul}(A)= Ext^{\bul}_{A\otimes A^{op}}(A,A)\,,
\end{equation}
where $A^{op}$ is the algebra $A$ with the opposite
multiplication and $A$ is naturally viewed as a left module
over $A\otimes A^{op}$\,.

Using the standard bar resolution for $A$ one shows that
Hochschild cohomology (\ref{Hoch}) is cohomology of the
following complex

\begin{equation}
\label{H-coch}
C^{m} (A)=Hom(A^{\otimes m}, A)\,,~(m\ge 1)\,, \qquad
C^{0}(A)=A\,
\end{equation}
with the differential given by
$$
\pa \Phi(a_1, \dots a_{m+1})=
a_1\cdot \Phi(a_2, \dots, a_{m+1}) -
\Phi(a_1\cdot a_2, a_3, \dots, a_{m+1})+ \dots
$$
\begin{equation}
\label{Hoch-diff}
+(-)^{m} \Phi(a_1,\dots, a_{m-1}, a_{m}\cdot a_{m+1})
+(-)^{m+1} \Phi(a_1,\dots, a_{m-1}, a_{m})\cdot a_{m+1}\,,
\end{equation}\\[-0.3cm]
$$
\pa \,:\, C^{m}(A) \mapsto C^{m+1}(A)\,.
$$
The vector space (\ref{H-coch}) is usually referred to
as the space of Hochschild cochains of the associative
algebra $A$\,.

Using the product in the algebra $A$ one can define
the associative cup product of Hochschild cochains.
The cup product of two homogeneous cochains
$\Phi_1\in C^{k_1}(A)$ and $\Phi_2\in C^{k_2}(A)$
is given by the formula
\begin{equation}
\label{cup-W}
\Phi_1 \cup \Phi_2(a_1, a_2,\dots, a_{k_1+k_2}) =
\Phi_1 (a_1, \dots, a_{k_1}) \cdot
\Phi_2 (a_{k_1+1},\dots, a_{k_1+k_2})\,,
\end{equation}
where $a_i \in A$\,.
It is not hard to show that the Hochschild differential
(\ref{Hoch-diff}) is a derivation of
the cup product (\ref{cup-W})
\begin{equation}
\label{cup-der}
\pa (\Phi_1\cup \Phi_2) = (\pa \Phi_1) \cup \Phi_2 +
(-)^{k_1}\Phi_1 \cup (\pa \Phi_2)
\end{equation}
and the induced product on the cohomology coincides
with the Yoneda product \cite{GM} in (\ref{Hoch})\,.

The space $C^{\bul}(A)$ can be also endowed with the so-called
Gerstenhaber bracket \cite{G} which is defined
between homogeneous elements $\Phi_1\in C^{k_1+1}(A)$ and
$\Phi_2\in C^{k_2+1}(A)$ as follows

$$[\Phi_1, \Phi_2]_G(a_0,\,\dots, a_{k_1+k_2})=$$
\begin{equation}
\label{Gerst-W}
\sum_{i=0}^{k_1}(-)^{ik_2}
\Phi_1(a_0,\,\dots , \Phi_2 (a_i,\,\dots,a_{i+k_2}),\, \dots, a_{k_1+k_2})
\end{equation}
$$
-(-)^{k_1k_2} (1 \leftrightarrow 2)\,, \qquad a_j \in A\,.
$$
Direct computation shows that (\ref{Gerst-W}) determines
a Lie (super)bracket on the space $C^{\bul}(A)[1]$ of the
Hochschild cochains with a shifted grading. One can observe that
the differential (\ref{Hoch-diff}) can be rewritten in terms
of the bracket (\ref{Gerst-W}) as follows

\begin{equation}
\label{pa}
\pa \Phi = (-)^{k+1}[\mu_0,\Phi]_G : C^k(A)\mapsto C^{k+1}(A)\,,
\end{equation}
where $\mu_0 \in C^2(A)$ is the multiplication
in the algebra $A$\,. This observation implies that
$\pa$ is a derivation of the Gerstenhaber bracket (\ref{Gerst-W})

\begin{equation}
\label{G-der-pa}
\pa[\Phi_1,\Phi_2]_G = [\pa \Phi_1,\Phi_2]_G +
(-)^{k_1-1} [\Phi_1, \pa \Phi_2]_G\,,
\qquad \Phi_i \in C^{k_i}(A)\,.
\end{equation}

\section{Appendix B. Equivariant resolution of the Weyl algebra.}
\renewcommand{\theequation}{B.\arabic{equation}}
In this section we propose a
$GL(2n,\bbC)$-equivariant homotopy
formula for the cohomological Hochschild complex
of the (formal) Weyl algebra\footnote{See paper \cite{Alev}
in which a similar computation has been performed for
the algebra of invariants of the ordinary (non-formal)
Weyl algebra acted upon by a finite group of automorphisms.}.

Let $\te^{ij}$ be a non-degenerate antisymmetric matrix of size
$2n\times 2n$. Then the vector space of the Weyl algebra $W$ is by definition
the vector space $\bbC[[y^1, \dots, y^{2n}]]((\h))$ of formal Laurent power series
in $\h$ whose coefficients are formal Taylor power series
in $y^1, \dots, y^{2n}$. The multiplication $\circ$ in $W$  is
given by

\begin{equation}
\label{circ-W}
a\circ b=\exp\left( \frac{\h}2\te ^{ij} \frac \partial {\partial
y^i}\frac \partial {\partial z^j}\right) a(\h, y)\,
b(\h,z)|_{z=y}\,, \qquad a,b\in W\,.
\end{equation}
One can observe that $W$ is naturally filtered with respect
to the degree of monomials $2[\h]+[y]$ where
$[\h]$ is a degree in $\h$ and $[y]$ is the total
degree in $y$'s

\begin{equation}
\label{filtr-W}
\dots \subset W^1  \subset W^0  \subset W^{-1}
\dots \subset W\,,
\quad
W^m = \{a\in W\,|\, a= \sum_{2k+p\ge m} \h^k a_{k; i_1\dots i_p}
y^{i_1} \dots y^{i_p} \}\,.
\end{equation}
This filtration defines the $2[\h]+[y]$-adic topology on the
algebra $W$ and the product (\ref{circ-W}) is continuous in
this topology.

Since $W$ is a topological algebra one should be careful
with the standard arguments of homological algebra.
In particular, the definition of Hochschild cohomology
for the algebra should be slightly modified.
By definition,
\begin{equation}
\label{Hoch-W}
HH^{k}(W)= H^{k} (C^{\bul}(W), \pa)\,, \qquad k\ge 0\,,
\end{equation}
where $C^{q}(W)$ is the vector space of continuous
$\bbC((\h))$-linear maps
$$
\Phi\in Hom_{\bbC((\h))}(W^{\totimes q}, W)\,,
$$
and $\totimes$ stands for the tensor product
over $\bbC((\h))$ completed in topology (\ref{filtr-W}).
Any such map can be uniquely represented in the
form of the following formal series

\begin{equation}
\label{operr-W}
\Phi =\sum_m \sum_{\al_1 \dots \al_q} \sum_{p=0}^{\infty}
\h^m \Phi^{\al_1\dots \al_q}_{m, i_1\dots i_p} y^{i_1}
\dots y^{i_p} \frac{\pa}{\pa y^{\al_1}}\otimes  \dots \otimes \frac{\pa}{\pa
y^{\al_q}}\,,
\end{equation}
where $\al$'s are multi-indices $\al=(j_1\dots j_l)$\,,
$$
\frac{\pa}{\pa y^{\al}}=\frac{\pa}{\pa y^{j_1}} \dots \frac{\pa}{\pa
y^{j_l}}\,,
$$
and the summation in $m$ is bounded below.

It is not hard to see that the Hochschild differential
$\pa$ (\ref{Hoch-diff}), the cup-product (\ref{cup-W}) and
the Gerstenhaber bracket (\ref{Gerst-W}) still make
sense for the complex $C^{\bul}(W)$.

Similarly, we replace the ``stupid'' bar resolution
by the topological bar resolution of $W$ as
a left $W\totimes W^{op}$-module

\begin{equation}
\label{bar}
B=\bigoplus_{m=0}^{\infty} B_m\,, \qquad
B_m= W^{\totimes (m+2)} = \bbC [[\by_1, \dots, \by_{m+2}]]((\h))\,,
\end{equation}
where $\by_i=(y^1_i, \dots, y^{2n}_i)$, $W^{op}$ denotes
the Weyl algebra with the opposite multiplication,
the differential $\mb\,:\,B_{m} \mapsto B_{m-1}$ is
given by the formula
$$
(\mb b)(\h,\by_1, \dots, \by_{m+1})=
$$
\begin{equation}
\label{mb-bar}
\sum_{k=1}^{m} (-)^{k-1} \exp(\frac{\h}2
\te^{ij}\frac{\pa}{\pa y_k^i} \frac{\pa}{\pa z^j})
b(\h, \by_1, \dots, \by_k, \bz, \by_{k+1}, \dots,
\by_{m+1})\Big|_{\bz=\by_k}\,,
\end{equation}
and the $W\totimes W^{op}$-module structure is defined by

$$
(a\cdot b)(\h,\by_1, \dots, \by_{m+2})= \exp(\frac{\h}2
\te^{ij}\frac{\pa}{\pa z^i} \frac{\pa}{\pa y_1^j})
a(\h,\bz) b(\h,\by_1, \dots, \by_{m+2})\Big|_{\bz=\by_1}\,,
$$
\begin{equation}
\label{WW-mod}
(b\cdot a)(\h,\by_1, \dots, \by_{m+2})= \exp(\frac{\h}2
\te^{ij}\frac{\pa}{\pa y_{m+2}^i} \frac{\pa}{\pa z^j})
b(\h,\by_1, \dots, \by_{m+2}) a(\h,\bz)\Big|_{\bz=\by_{m+2}}\,,
\end{equation}
$$
a\in W,\qquad b \in \bbC[[\by_1, \dots, \by_{m+2}]]((\h))\,.
$$
By $h_B$ we denote the corresponding homotopy operator
$h_B \,:\, B_{m-1} \mapsto B_m$
\begin{equation}
\label{bar-h}
(h_B b) (\h,\by_1 ,\dots, \by_{m+2}) =
b(\h,\by_2 ,\dots, \by_{m+2})
\end{equation}
in the augmented bar complex

\begin{equation}
\label{aug-bar}
\dots \stackrel{\mb}{\rightarrow} B_2
  \stackrel{\mb}{\rightarrow}  B_1
   \stackrel{\mb}{\rightarrow} B_0
   \stackrel{\circ}{\rightarrow} W\,.
\end{equation}

Obviously, we have the following canonical isomorphism
of complexes

\begin{equation}
\label{Hoch-def}
C^{\bul}(W)= Hom^{c}_{W\totimes W^{op}} (B_{\bul}, W)\,,
\end{equation}
where by $Hom^c$ we denote the vector space of all
homomorphisms continuous in topology (\ref{filtr-W}).

An appropriate analogue of the Koszul resolution
of the Weyl algebra $W$ is the polynomial algebra

\begin{equation}
\label{Koszul}
K= W\totimes W^{op} [C^1, \dots, C^{2n}]
\end{equation}
in $2n$ {\it anticommuting} variables $C^1, \dots, C^{2n}$
with coefficients in the tensor product
$W\totimes W^{op}$ over $\bbC((\h))$ completed
in topology (\ref{filtr-W}). The differential in $K$ is given by
the formula
\begin{equation}
\label{Koszul-d}
\md a = (y^i_1- y^i_2) \frac{\vec{\pa}}{\pa C^i} a+
\frac{\h}{2} \te^{ij} \frac{\vec{\pa}}{\pa C^i}
(\frac{\pa}{\pa y_1^i}+\frac{\pa}{\pa y_2^i}) a\,,
\qquad
a\in K\,,
\end{equation}
where the arrow over $\pa$ means that we use
left derivatives with respect to anticommuting
variables $C^i$\,.

The $m$-th term of the complex (\ref{Koszul})
$$
K_m=\wedge^{m}(\bbC^{2n}) \otimes W\totimes W^{op}
$$
is naturally a left module over $W\totimes W^{op}$
and the differential (\ref{Koszul-d})
$\md\,:\, K_m \mapsto K_{m-1}$ is compatible with
the $W\totimes W^{op}$-module structure.

The homotopy operator $h_K$ of the augmented Koszul complex

\begin{equation}
\label{Koszul-aug}
\dots \stackrel{\md}{\rightarrow} \wedge^m(\bbC^{2n})\otimes W\totimes W^{op}
\stackrel{\md}{\rightarrow}
 \dots \stackrel{\md}{\rightarrow} \bbC^{2n}\otimes W\totimes W^{op}
   \stackrel{\md}{\rightarrow} W\totimes W^{op}  \stackrel{\circ}{\rightarrow} W
\end{equation}
looks as follows
\begin{equation}
\label{Koszul-h}
\begin{array}{c}
\displaystyle
h_K(a) = C^k \int_0^1 dt (\cD_{-t}\, \cD\,
\frac{\pa}{\pa y_1^k}a)((\h, \by_2+ t(\by_1 -\by_2), \by_2,
tC))\,, \\[0.5cm]
\displaystyle
\cD= \exp\left( \frac{\h}{2}\te^{ij}\frac{\pa}{\pa y_1^i}\frac{\pa}{\pa y_2^j}
\right)\,,\\[0.5cm]
\displaystyle
\cD_{-t}=
\exp\left( -\frac{\h t}{2}\te^{ij}\frac{\pa}{\pa y_1^i}\frac{\pa}{\pa
y_2^j}\right)\,.
\end{array}
\end{equation}

Using homotopy operators (\ref{bar-h}), (\ref{Koszul-h})
we can inductively construct continuous quasi-isomorphisms
between the topological bar resolution (\ref{bar}) and
the Koszul resolution (\ref{Koszul}). Namely,
\begin{pred}\label{K-bar-K}
The complexes of left $W\totimes W^{op}$-modules
(\ref{bar}) and (\ref{Koszul}) are quasi-isomorphic.
A quasi-isomorphism $\la$ from the Koszul resolution
(\ref{Koszul})
to the topological bar resolution
(\ref{bar}) is determined inductively by setting
$$
\la\Big|_{W\totimes W^{op}}= id \,:\, K_0 \mapsto B_0\,,
$$
and
$$
\la (a) = h_B (\la (\md a))
$$
for any topological generator $a$ of the left
$W\totimes W^{op}$-module $K_m$ $(m>0)$\,.
A quasi-isomorphism $\nu$ from the topological
bar resolution (\ref{bar}) to the Koszul resolution
(\ref{Koszul}) is determined inductively by setting
$$
\nu\Big|_{W\totimes W^{op}}= id \,:\, B_0 \mapsto K_0\,,
$$
and
$$
\nu (b) = h_K (\nu (\mb b))
$$
for any topological generator $b$ of the left
$W\totimes W^{op}$-module $B_m$ $(m>0)$\,.
\end{pred}
{\bf Proof.} We define inductively
a morphism $\rho\,:\, B\mapsto B[1]$ of
left $W\totimes W^{op}$-modules
by setting

\begin{equation}
\label{rho-def}
\rho\Big|_{W\totimes W^{op}}= 0 \,:\, B_0 \mapsto B_1\,,
\qquad
\rho (b) = h_B (id-\la \nu)(b) - h_B \rho(\mb(b))
\end{equation}
for any topological generator $b$ of the left
$W\totimes W^{op}$-module $B_m$ $(m>0)$\,.
Direct computation shows that $\rho$ is
a homotopy operator between $id$ and
$\la\, \nu$

\begin{equation}
\label{rho-prop}
b- \la (\nu(b)) = \mb (\rho(b)) + \rho(\mb(b))\,.
\end{equation}
Similarly, one constructs a homotopy operator
between $id$ and $\nu\, \la$ on $K$ $\Box$.

The above proposition immediately implies that
the maps
\begin{equation}
\label{la-hat}
(\hat\la a)(\ka) = a (\la \ka) \,:\,
C^{\bul}(W) \mapsto Hom^c_{W\totimes W^{op}}(K_{\bul},W)\,,
\qquad \ka \in K_{\bul}\,,~ a\in C^{\bul}(W)
\end{equation}

\begin{equation}
\label{nu-hat}
(\hat\nu f)(b) = f(\nu b) \,:\,
Hom^c_{W\totimes W^{op}}(K_{\bul},W)
\mapsto C^{\bul}(W)\,,
\end{equation}
$$
b \in B_{\bul}\,,\qquad  f\in Hom^c_{W\totimes W^{op}}(K_{\bul},W)
$$
are quasi-isomorphisms of
the complexes $C^{\bul}(W)$ and
$Hom^c_{W\totimes W^{op}}(K_{\bul},W)$ and the map

\begin{equation}
\label{rho-hat}
(\hat\rho a)(b) = a (\rho b) \,:\,
C^{\bul}(W) \mapsto C^{\bul}(W)\,,
\qquad  a \in C^{\bul}(W)\,,~ b \in B_{\bul}
\end{equation}
is a homotopy operator between the $id$ and $\hat\nu\,\hat\la$

\begin{equation}
\label{rho-hat-prop}
a - \hat\nu (\hat\la\, a) =
(\pa \hat\rho + \hat\rho \pa) a\,,
\qquad \forall~ a \in C^{\bul}(W)\,.
\end{equation}
In the above equations we use the identification
of $C^{\bul}(W)$ and $Hom^c_{W\totimes W^{op}}(B_{\bul}, W)$
(\ref{Hoch-def}).

The complex $Hom^c_{W\totimes W^{op}}(K_{\bul},W)$
is obviously identified with the polynomial
algebra $W[\psi_1, \dots, \psi_{2n}]$
in $2n$ {\it anticommuting} variables
$\psi_1, \dots, \psi_{2n}$ with coefficients in
the Weyl algebra $W$, namely

\begin{equation}
\label{Hom-W}
Hom^c_{W\totimes W^{op}}(K_{q},W) = W[\psi_1, \dots, \psi_{2n}]_{(q)},
\end{equation}
where $W[\psi_1, \dots, \psi_{2n}]_{(q)}$ is the
space of polynomials in $\psi$'s of degree $q$.
One can easily compute
the differential $\pa_H$ in $ W[\psi_1, \dots, \psi_{2n}]$

\begin{equation}
\label{pa-H}
\pa_H a = \h \psi_i \om^{ij}\frac{\pa}{\pa y^j} a\,,
\qquad
a\in  W[\psi_1, \dots, \psi_{2n}],
\end{equation}
and guess the corresponding partial homotopy operator

\begin{equation}
\label{homot-H}
\cH a = \frac1{\h}
\int_0^1 d t\,
y^i \om_{ij} (\frac{\vec{\pa}}{\pa \psi_j} a)(\h, t y, t\psi)\,,
\qquad
a\in  W[\psi_1, \dots, \psi_{2n}]
\end{equation}
satisfying the following identity

\begin{equation}
\label{homot-prop}
a = a\Big|_{y=\psi=0} + (\pa_H \cH + \cH \pa_H) a  \,,
\qquad
\forall~a\in  W[\psi_1, \dots, \psi_{2n}]\,.
\end{equation}
Thus we conclude that

\begin{equation}
\label{su}
HH^{q}(W)=\begin{cases}
\bbC((\h))\,,~ q=0\,,\cr
0\,,~ {\rm otherwise}\,.
\end{cases}
\end{equation}
For the case $q=0$ we also have an obvious
algebraic formulation of the above assertion

\begin{equation}
\label{0-cocyc}
Z^0(C^{\bul}(W), \pa)= Z(W)
= \bbC((\h))\,,
\end{equation}
where $Z^0(C^{\bul}(W), \pa)$ denotes
the space of the zeroth cocycles of the
complex $C^{\bul}(W)$ and $Z(W)$ stands for
the center of the Weyl algebra $W$.

Together with the expected result (\ref{su})
the above considerations enable us to
get a homotopy formula for the Hochschild
cohomological complex of the Weyl algebra.
Using proposition \ref{K-bar-K} and
equations (\ref{rho-hat-prop}), (\ref{homot-prop})
we derive that
\begin{pred}
\label{homot-Hoch}
The following operator
\begin{equation}
\label{xxx}
\chi (a) = \hat\nu (\cH (\hat\la(a))) + \hat\rho(a)\,,
\qquad a\in C^{\ge 1}(W)\,
\end{equation}
satisfies the equation
\begin{equation}
\label{xxx-prop}
a= (\pa \chi + \chi \pa) a
\end{equation}
for any $a\in C^{\ge 1}(W)$\,. $\Box$
\end{pred}
An explicit expression of the homotopy operator
(\ref{xxx}) is rather complicated.
See, for example, paper \cite{H}
in which an explicit homotopy formula is given for
the Hochschild complexes of the commutative algebra
of polynomials in $n$ variables.

One can observe that any element $g\in GL(2n,\bbC)$
defines an isomorphism from the Weyl algebra $W_{\te}$
associated to the (antisymmetric, non-degenerate)
matrix $\te$ to the Weyl algebra
$W_{\te'}$ associated to the matrix $\te'= g_{\ast}\te$.
For our purposes we need to keep track of the
equivariance properties with respect to such
isomorphisms. To do this we introduce a groupoid of Weyl algebras.
The objects of this groupoid are Weyl algebras
$W_{\te}$ associated to all possible antisymmetric
non-degenerate matrices $\te$ and the morphisms
are defined by elements of $GL(2n,\bbC)$\,.
One can repeat all the arguments in this appendix
by replacing a single Weyl algebra by this groupoid.
Thus we get the following
\begin{pred} \label{equiv}
The homotopy operator (\ref{xxx}) is
$GL(2n,\bbC)$-equivariant. Namely, if
$\chi_{\te}$ denotes the homotopy operator
(\ref{xxx}) in the complex $C^{\bul}(W_{\te})$
then for any $g\in GL(2n, \bbC)$ the diagram
\begin{equation}
\label{diag-equiv}
\begin{array}{ccc}
C^{\bul}(W_{\te}) &\stackrel{\chi_{\te}}{\rightarrow} & C^{\bul-1}(W_{\te})\\[0.3cm]
\downarrow^{g_{\ast}}  & ~  &     \downarrow^{g_{\ast}} \\[0.3cm]
C^{\bul}(W_{g_{\ast}\te})  &\stackrel{\chi_{g_{\ast}\te}}{\rightarrow} &
 C^{\bul-1}(W_{g_{\ast}\te})
\end{array}
\end{equation}
commutes. $\Box$
\end{pred}

\end{document}